\documentclass[11pt]{article}
\usepackage{amsfonts}
\usepackage{amsmath}
\usepackage{amssymb}
\usepackage{graphicx}
\usepackage[dvipsnames,usenames]{color}
\usepackage[pdfstartview=FitH ]{hyperref}
\usepackage[hyperpageref]{backref}
\begin{document}

\baselineskip 18pt
\def\o{\over}
\def\e{\varepsilon}
%------------------------------------------------------------------
\title{\Large\bf On\ \ sumsets\ \ in\ \ ${\Bbb F}_2^n$}
\author{Chaohua\ \  Jia}
\date{}
\maketitle {\small \noindent {\bf Abstract.} Let ${\Bbb F}_2$ be the
finite field of two elements, ${\Bbb F}_2^n$ be the vector space of
dimension $n$ over ${\Bbb F}_2$. For sets $A,\,B\subseteq{\Bbb
F}_2^n$, their sumset is defined as the set of all pairwise sums
$a+b$ with $a\in A,\,b\in B$.

Ben Green and Terence Tao proved that, let $K\geq 1$, if
$A,\,B\subseteq{\Bbb F}_2^n$ and $|A+B|\leq K|A|^{1\o 2}|B|^{1\o
2}$, then there exists a subspace $H\subseteq{\Bbb F}_2^n$ with
$$
|H|\gg\exp(-O(\sqrt{K}\log K))|A|
$$
and $x,\,y\in{\Bbb F}_2^n$ such
that
$$
|A\cap(x+H)|^{1\o 2}|B\cap(y+H)|^{1\o 2}\geq{1\o 2K}|H|.
$$

In this note, we shall use the method of Green and Tao with some
modification to prove that if
$$
|H|\gg\exp(-O(\sqrt{K}))|A|,
$$
then the above conclusion still holds true.
}

%------------------------------------------------------------------
\vskip.5in
\noindent{\bf 1. Introduction}

Let ${\Bbb F}_2$ be the finite field of two elements, ${\Bbb F}_2^n$
be the vector space of dimension $n$ over ${\Bbb F}_2$. For sets
$A,\,B\subseteq{\Bbb F}_2^n$, their sumset $A+B$ is defined as
$$
A+B:=\{a+b:\ a\in A,\,b\in B\}.
$$

In 1999, Ruzsa[4] proved the following theorem.

{\bf Theorem 1}(Ruzsa). Let $K\geq 1$ be an integer, and suppose
that set $A\subseteq{\Bbb F}_2^n$ with $|A+A|\leq K|A|$. Then $A$ is
contained in a subspace $H\subseteq{\Bbb F}_2^n$ with $|H|\leq
F(K)|A|$, where $F(K)=K^22^{K^4}$.

This result was improved by Sanders[5] to $F(K)=2^{O(K^{3\o 2}\log
K)}$ in 2008 and then improved by Green and Tao[2] to
$F(K)=2^{2K+O(\sqrt{K}\log K)}$ in 2009. The bound
$F(K)=2^{2K+O(\sqrt{K}\log K)}$ is almost best possible.

If we do not require that the subspace $H$ contains the set $A$
completely but contains a part of $A$, then related bounds can be
further improved. The following theorem was given in [1] and some
explanations on it could be found in the introduction of [3].

{\bf Theorem 2}. Suppose that $K\geq 1$ and that $A\subseteq{\Bbb
F}_2^n$ with $|A+A|\leq K|A|$. Then there is a subspace
$H\subseteq{\Bbb F}_2^n$ with $|H|\ll K^{O(1)}|A|$ such that
$$
|A\cap H|\gg\exp(-K^{O(1)})|A|.
$$

If we permit to replace the subspace $H$ by translates of it, then
better bounds could be obtained. In 2009, Green and Tao[3] obtained
the following result.

{\bf Theorem 3}(Green-Tao). Let $K\geq 1$, if $A,\,B\subseteq{\Bbb
F}_2^n$ and $|A+B|\leq K|A|^{1\o 2}|B|^{1\o 2}$, then there exists a
subspace $H\subseteq{\Bbb F}_2^n$ with
$$
|H|\gg\exp(-O(\sqrt{K}\log K))|A|
$$
and $x,\,y\in{\Bbb F}_2^n$ such that
$$
|A\cap(x+H)|^{1\o 2}|B\cap(y+H)|^{1\o 2}\geq{1\o 2K}|H|.
$$

In this note, we shall use the method of Green and Tao with some
modification to prove the following theorem.

{\bf Theorem 4}. Let $K\geq 1$, if $A,\,B\subseteq{\Bbb F}_2^n$ and
$|A+B|\leq K|A|^{1\o 2}|B|^{1\o 2}$, then there exists a subspace
$H\subseteq{\Bbb F}_2^n$ with
$$
|H|\gg\exp(-O(\sqrt{K}))|A|
$$
and $x,\,y\in{\Bbb F}_2^n$ such that
$$
|A\cap(x+H)|^{1\o 2}|B\cap(y+H)|^{1\o 2}\geq{1\o 2K}|H|.
$$

\vskip.3in
\noindent{\bf 2. Definitions}

In this section we shall introduce some definitions given in [3].

{\bf Definition 1}(normalized energy). For non-empty sets
$A_1,\,A_2,\,A_3,\,A_4\subseteq {\Bbb F}_2^n$, define the normalized
energy \begin{align*} &\omega(A_1,\,A_2,\,A_3,\,A_4):={1\o
(|A_1||A_2| |A_3||A_4|)^{3\o 4}}\,|\{(a_1,\,a_2,\,a_3,\,a_4)\\
&\qquad\in A_1\times A_2\times A_3\times A_4:\ a_1+a_2+a_3+a_4=0\}|.
\end{align*}

It was shown in [3] that
\begin{align}
0\leq\omega(A_1,\,A_2,\,A_3,\,A_4)\leq 1.
\end{align}

{\bf Definition 2}(Fourier transform). For $f:\,{\Bbb
F}_2^n\longrightarrow {\Bbb R}$, define the Fourier transform
$\hat{f}:\,{\Bbb F}_2^n\longrightarrow {\Bbb R}$ by
$$
\hat{f}(\xi):={1\o 2^n}\sum_{x\in{\Bbb F}_2^n}f(x)(-1)^{\xi\cdot x},
$$
where
$$
\xi\cdot x=(\xi_1,\,\cdots,\,\xi_n)\cdot(x_1,\,\cdots,\,x_n)=\xi_1
x_1+\cdots +\xi_n x_n.
$$

{\bf Definition 3}(spectrum). If $A\subseteq {\Bbb F}_2^n$ is
non-empty and $0<\alpha \leq 1$, define the $\alpha-$spectrum
$$
{\rm Spec}_\alpha(A):=\{\xi\in{\Bbb F}_2^n:\ |\hat{\bf
1}_A(\xi)|\geq \alpha\,{|A|\o 2^n}\},
$$
where ${\bf 1}_A(x)$ is the indicator function of set $A$.

{\bf Definition 4}(coherently flat quadruples). Suppose that
$A_1,\,A_2,\,A_3,\,A_4\subseteq {\Bbb F}_2^n$ are non-empty and
$\delta\in(0,\,{1\o 2})$ is a small parameter. If for each $\xi\in
{\Bbb F}_2^n$, one of the following conditions is satisfied:

1) $\xi\in{\rm Spec}_{9\o 10}(A_i)$ for all $i=1,\,2,\,3,\,4$;

2) $\xi\not\in{\rm Spec}_\delta(A_i)$ for all $i=1,\,2,\,3,\,4$,

\noindent we say that the quadruple $(A_1,\,A_2,\,A_3,\,A_4)$ is
coherently $\delta-$ flat.

\vskip.3in
\noindent{\bf 3. The proof of Theorem 4}

{\bf Lemma 1}. Let $J\geq 1$. Suppose that $(A_1,\,A_2,\,A_3,\,A_4)$
is a coherently ${1\o \sqrt{2J}}-$flat quadruple, the normalized
energy of which satisfies
$$
\omega(A_1,\,A_2,\,A_3,\,A_4)\geq {1\o J}.
$$
Then there is a subspace $H\subseteq{\Bbb F}_2^n$ with
$x_1,\,x_2,\,x_3,\,x_4 \in{\Bbb F}_2^n$ such that
\begin{align}
H\geq{4\o 5}\,(|A_1| |A_2| |A_3| |A_4|)^{1\o 4}
\end{align}
and
\begin{align}
\prod_{i=1}^4|A_i\cap(x_i+H)|^{1\o 4}\geq {1\o 2J}\,|H|.
\end{align}

This is Proposition 2.4 in [3].

Let
$$
{\rm Dbl}(A,\,B):={|A+B|\o |A|^{1\o 2}|B|^{1\o 2}}.
$$
Since
$$
|A+B|\geq\max(|A|,\,|B|),
$$
we have
\begin{align}
{\rm Dbl}(A,\,B)\geq 1.
\end{align}

{\bf Lemma 2}. Suppose that $A,\,B\subseteq{\Bbb F}_2^n$ are
non-empty and that for $J\geq 1$,
$$
{\rm Dbl}(A,\,B)\leq J.
$$
If $(A,\,B,\,A,\,B)$ is not coherently ${1\o \sqrt{2J}}-$flat, then
there are $A' \subseteq A,\,B'\subseteq B$ such that
\begin{align}
|A'|\geq {1\o 20}\,|A|,\qquad\qquad |B'|\geq {1\o 20}\,|B|
\end{align}
and
\begin{align}
{\rm Dbl}(A',\,B')\leq{J\o 1+{1\o {100\sqrt{J}}}}.
\end{align}

{\bf Proof}. By the supposition, there is $\xi\in{\Bbb F}_2^n$ such
that
\begin{align}
\xi\not\in{\rm Spec}_{9\o 10}(A)\cap{\rm Spec}_{9\o 10}(B)
\end{align}
and
\begin{align}
\xi\in{\rm Spec}_{1\o \sqrt{2J}}(A)\cup{\rm Spec}_{1\o
\sqrt{2J}}(B).
\end{align}
By (7), $\xi\ne 0$. Write
\begin{align*}
&A_0:=\{x\in A:\ x\cdot\xi=0\},\qquad\quad  A_1:=\{x\in A:\ x\cdot
\xi=1\},\\
&B_0:=\{x\in B:\ x\cdot\xi=0\},\qquad\quad  B_1:=\{x\in B:\ x\cdot
\xi=1\}.
\end{align*}

If $|A_0|\geq{1\o 2}|A|$, we write $\alpha:={|A_0|\o |A|}$.
Otherwise, $|A_0| <{1\o 2}|A|\Rightarrow |A_1|=|A|-|A_0|\geq
|A|-{1\o 2}|A|={1\o 2}|A|$. Then we write $\alpha:={|A_1|\o |A|}$.
Without loss of generality, we can suppose that $|A_0|\geq{1\o
2}|A|$ and write
$$
\alpha:={|A_0|\o |A|}.
$$
Similarly, we can also suppose that $|B_0|\geq{1\o 2}|B|$ and write
$$
\beta:={|B_0|\o |B|}.
$$
We have
\begin{align}
\alpha\geq {1\o 2},\qquad\qquad {\beta}\geq{1\o 2}.
\end{align}

By
\begin{align*}
|\hat{\bf 1}_A(\xi)|&=\Bigl|{1\o 2^n}\sum_{x\in
A}(-1)^{x\cdot\xi}\Bigr|\\
&=\Bigl|{1\o 2^n}\Bigl(\sum_{x\in A_0}(-1)^{x\cdot\xi}+\sum_{x\in
A_1}(-1)^{x\cdot\xi}\Bigr)\Bigr|\\
&=\Bigl|{1\o 2^n}\,(|A_0|-|A_1|)\Bigr|\\
&=\Bigl|{1\o 2^n}\,(2|A_0|-|A|)\Bigr|\\
&=(2\alpha-1)\cdot{|A|\o 2^n}
\end{align*}
and
$$
|\hat{\bf 1}_B(\xi)|=(2\beta-1)\cdot{|B|\o 2^n},
$$
we know that the condition (7) is equivalent to
\begin{align}
2\alpha-1<{9\o 10}\qquad\quad{\rm or}\qquad\quad 2\beta-1<{9\o 10},
\end{align}
and the condition (8) is equivalent to
\begin{align}
2\alpha-1\geq{1\o \sqrt{2J}}\qquad\quad{\rm or}\qquad\quad 2
\beta-1\geq {1\o \sqrt{2J}}.
\end{align}

Without loss of generality, we suppose that
\begin{align}
\beta\geq\alpha
\end{align}
and consider
$$
|B_0+A_0|+|B_0+A_1|.
$$
If $\beta<\alpha$, we shall consider $|A_0+B_0|+|A_0+B_1|$.

It is easy to see that sets $B_0+A_0$ and $B_0+A_1$ are disjoint.
Hence,
$$
|B_0+A_0|+|B_0+A_1|\leq |B+A|\leq J|B|^{1\o 2}|A|^{1\o 2}
$$
or
\begin{align}
{|B_0+A_0|\o |B_0|^{1\o 2}|A_0|^{1\o 2}}\cdot{|B_0|^{1\o 2}|A_0|^
{1\o 2}\o |B|^{1\o 2}|A|^{1\o 2}}+{|B_0+A_1|\o |B_0|^{1\o 2}|A_1|
^{1\o 2}}\cdot{|B_0|^{1\o 2}|A_1|^{1\o 2}\o |B|^{1\o 2}|A|^{1\o 2}}
\leq J.
\end{align}
Let
$$
\Psi:=\min\Bigl({|B_0+A_0|\o |B_0|^{1\o 2}|A_0|^{1\o
2}},\,{|B_0+A_1| \o |B_0|^{1\o 2}|A_1|^{1\o 2}}\Bigr).
$$
It follows from (13) that
\begin{align}
\Psi(\beta^{1\o 2}\alpha^{1\o 2}+\beta^{1\o 2}(1-\alpha)^{1\o 2})
\leq J.
\end{align}

Under the supposition (12),  the condition (10) is equivalent to
\begin{align}
\alpha<{19\o 20},
\end{align}
and the condition (11) is equivalent to
\begin{align}
\beta\geq {1\o 2}+{1\o 2\sqrt{2J}}.
\end{align}
We shall discuss in the following two cases.

Case 1. ${1\o 2}+{1\o 2\sqrt{2J}}\leq\alpha<{19\o 20}$.

By (12),
$$
\beta^{1\o 2}\alpha^{1\o 2}+\beta^{1\o 2}(1-\alpha)^{1\o 2}
\geq\alpha+\alpha^{1\o 2}(1-\alpha)^{1\o 2}.
$$
The discussion in [3] yields that
\begin{align*}
&\
\,\alpha+\alpha^{1\o 2}(1-\alpha)^{1\o 2}\geq\alpha+2\alpha
(1-\alpha)\\
&=1+(2\alpha-1)(1-\alpha)\geq 1+{1\o 20\sqrt{2J}}.
\end{align*}
Hence,
$$
\beta^{1\o 2}\alpha^{1\o 2}+\beta^{1\o 2}(1-\alpha)^{1\o 2} \geq
1+{1\o 20\sqrt{2J}}.
$$

Case 2. ${1\o 2}\leq\alpha<{1\o 2}+{1\o 2\sqrt{2J}}$.

It follows from (16) that
$$
\beta^{1\o 2}\alpha^{1\o 2}+\beta^{1\o 2}(1-\alpha)^{1\o 2}
\geq\Bigl({1\o 2}+{1\o 2\sqrt{2J}}\Bigr)^{1\o 2}(\alpha^{1\o 2}+
(1-\alpha)^{1\o 2}).
$$
Let
$$
f(\alpha)=\alpha^{1\o 2}+(1-\alpha)^{1\o 2}.
$$
Since
$$
f'(\alpha)={1\o 2\sqrt{\alpha}}-{1\o 2\sqrt{1-\alpha}}\leq 0,
$$
the function $f(\alpha)$ is decreasing monotonically. Thus,
$$
\alpha^{1\o 2}+(1-\alpha)^{1\o 2}\geq\Bigl({1\o 2}+{1\o 2\sqrt{2J}}
\Bigr)^{1\o 2}+\Bigl(1-\Bigl({1\o 2}+{1\o 2\sqrt{2J}}\Bigr)\Bigr)
^{1\o 2}.
$$
Hence,
\begin{align*}
&\ \,\beta^{1\o 2}\alpha^{1\o 2}+\beta^{1\o 2}(1-\alpha)^{1\o 2}\\
&\geq\Bigl({1\o 2}+{1\o 2\sqrt{2J}}\Bigr)^{1\o 2}\Bigl(\Bigl({1\o
2}+ {1\o 2\sqrt{2J}}\Bigr)^{1\o 2}+\Bigl(1-\Bigl({1\o 2}+{1\o
2\sqrt{2J}}\Bigr)\Bigr)^{1\o 2}\Bigr)\\
&=\Bigl({1\o 2}+ {1\o 2\sqrt{2J}}\Bigr)+\Bigl({1\o 2}+{1\o
2\sqrt{2J}}\Bigr)^{1\o 2}\Bigl(1-\Bigl({1\o 2}+{1\o
2\sqrt{2J}}\Bigr)\Bigr)^{1\o 2}
\end{align*}
which is the value of function $\alpha+\alpha^{1\o 2}(1-\alpha)^{1\o
2}$ at $\alpha={1\o 2} +{1\o 2\sqrt{2J}}$. By the discussion in Case
1, we have
$$
\beta^{1\o 2}\alpha^{1\o 2}+\beta^{1\o 2}(1-\alpha)^{1\o 2} \geq
1+{1\o 20\sqrt{2J}}.
$$

Combining the above two cases, we get
$$
\Psi\leq{J\o 1+{1\o 20\sqrt{2J}}}\leq{J\o 1+{1\o 100\sqrt{J}}}.
$$
Take $B'=B_0,\,A'=A_0$ or $A_1$ such that
$$
\Psi={\rm Dbl}(A',\,B').
$$
Then
$$
{\rm Dbl}(A',\,B')\leq{J\o 1+{1\o 100\sqrt{J}}}.
$$
Since
$$
|A_0|\geq{1\o 2}|A|,\qquad\quad |A_1|=|A|-|A_0|\geq |A|-{19\o
20}|A|={1\o 20}|A|,
$$
we have
$$
|A'|\geq{1\o 20} |A|.
$$
We also have
$$
|B'|\geq{1\o 20} |B|.
$$

So far the proof of Lemma 2 is finished.

{\bf Lemma 3}. Suppose that $A,\,B\subseteq{\Bbb F}_2^n$ are
non-empty and that for $K\geq 1$,
$$
{\rm Dbl}(A,\,B)\leq K.
$$
Then there are $A' \subseteq A,\,B'\subseteq B$ with
\begin{align}
|A'|\gg\exp(-O(\sqrt{K}))|A|,\qquad |B'|\gg\exp(-O(\sqrt{K}))|B|
\end{align}
such that for some $J(1\leq J\leq K)$,
\begin{align}
{\rm Dbl}(A',\,B')\leq J
\end{align}
and $(A',\,B',\,A',\,B')$ is coherently ${1\o \sqrt{2J}}-$flat.

{\bf Proof}. Take $K_1=K$. If $(A,\,B,\,A,\,B)$ is coherently ${1\o
\sqrt{2K}}-$flat, then the conclusion holds true.

If $(A,\,B,\,A,\,B)$ is not coherently ${1\o \sqrt{2K}}-$flat, Lemma
2 produces that there are $A''\subseteq A,\,B''\subseteq B$ with
$$
|A''|\geq {1\o 20}|A|,\qquad\quad |B''|\geq {1\o 20}|B|
$$
such that
$$
{\rm Dbl}(A'',\,B'')\leq{K_1\o 1+{1\o 100\sqrt{K_1}}}.
$$
Then take
$$
K_2={K_1\o 1+{1\o 100\sqrt{K_1}}},
$$
and for $A'',\,B''$ and $K_2$, repeat the above process.

Since ${\rm Dbl}\geq 1$, this process has to stop after finite
steps. We get a sequence $K_1=K,\,K_2,\,\cdots,\,K_m=J$ with
$$
K_{i+1}={K_i\o 1+{1\o 100\sqrt{K_i}}},\qquad\qquad
i=1,\,2,\,\cdots,\,m-1
$$
and $A'\subseteq A,\,B'\subseteq B$ with
$$
|A'|\gg{1\o (20)^m}|A|,\qquad\quad |B'|\gg{1\o (20)^m}|B|
$$
such that
$$
{\rm Dbl}(A',\,B')\leq J
$$
and $(A',\,B',\,A',\,B')$ is coherently ${1\o \sqrt{2J}}-$flat.

We distribute $K_i$ into intervals
$$
\Bigl({K\o e^{r+1}},\,{K\o e^r}\Bigr],\,\Bigl({K\o e^r},\,{K\o
e^{r-1}}\Bigr],\,\cdots,\,\Bigl({K\o e^2},\,{K\o e}\Bigr],\,
\Bigl({K\o e},\,K\Bigr],\quad r=[\log K].
$$
For the given interval $({K\o e^{s+1}},\,{K\o e^s}] (0\leq s\leq
r)$, if $K_l$ and $K_{l+j}(j\geq 1)\in ({K\o e^{s+1}},\,{K\o e^s}]$,
we have
\begin{align*}
{K\o e^{s+1}}&\leq K_{l+j}={K_{l+j-1}\o 1+{1\o
100\sqrt{K_{l+j-1}}}}\leq{K_{l+j-1}\o 1+{1\o 100\sqrt{K\o
e^s}}}\leq\cdots\\
&\leq{K_l\o \Bigl(1+{1\o 100\sqrt{K\o e^s}}\Bigr)^j}\leq{K\o
e^s}\cdot{1\o \Bigl(1+{1\o 100\sqrt{K\o e^s}}\Bigr)^j}.
\end{align*}
Thus
\begin{align*}
&\ \Bigl(1+{1\o 100\sqrt{K\o e^s}}\Bigr)^j\leq e,\\
&j\cdot{1\o {\sqrt{K\o e^s}}}\ll j\log\Bigl(1+{1\o 100\sqrt{K\o
e^s}}\Bigr)\leq 1,\\
&\qquad\quad j\ll\sqrt{K\o e^s}.
\end{align*}
Hence, the number of $K_i$ dropping into the interval $({K\o
e^{s+1}},\,{K\o e^s}]$ is $\ll\sqrt{K\o e^s}$. For the total number
of $K_i$, we have
\begin{align*}
m&\ll\sqrt{K}+\sqrt{K\o e}+\sqrt{K\o e^2}+\cdots+\sqrt{K\o e^r}\\
&\leq\sqrt{K}\Bigl(1+{1\o \sqrt{e}}+{1\o (\sqrt{e})^2}+{1\o
(\sqrt{e})^3}+\cdots\Bigr)\\
&\ll\sqrt{K}.
\end{align*}
Therefore
$$
|A'|\gg\exp(-O(\sqrt{K}))|A|,\qquad |B'|\gg\exp(-O(\sqrt{K}))|B|.
$$

So far the proof of Lemma 3 is finished.

{\bf The proof of Theorem 4}. We take $A',\,B'$ in Lemma 3 with
required properties. It is shown in [3] that
$$
\omega(A',\,B',\,A',\,B')\geq{1\o {\rm Dbl}(A',\,B')}\geq{1\o J}.
$$
Lemma 1 claims that there is a subspace $H\subseteq{\Bbb F}_2^n$
with $x_1,\,x_2,\,x_3,\,x_4 \in{\Bbb F}_2^n$ such that
$$
H\geq{4\o 5}\,|A'|^{1\o 2}|B'|^{1\o 2}\gg \exp(-O(\sqrt{K}))|A|^
{1\o 2}|B|^{1\o 2}
$$
and
\begin{align*}
&\ \,|A\cap(x_1+H)|^{1\o 4}|B\cap(x_2+H)|^{1\o 4}|A\cap(x_3+H)|^{1\o
4}|B\cap(x_4+H)|^{1\o 4}\\
&\geq {1\o 2J}\,|H|\geq {1\o 2K}\,|H|.
\end{align*}

Since
$$
|A|\leq |A+B|\leq K|A|^{1\o 2}|B|^{1\o 2},
$$
we have
$$
K^{-2}|A|\leq |B|,
$$
hence
$$
H\gg\exp(-O(\sqrt{K}))|A|.
$$

Without loss of generality, we can suppose that
$$
|A\cap(x_1+H)|\geq |A\cap(x_3+H)|,\quad |B\cap(x_2+H)|\geq |B\cap
(x_4+H)|,
$$
hence
$$
|A\cap(x_1+H)|^{1\o 2}|B\cap(x_2+H)|^{1\o 2}\geq {1\o 2K}\,|H|.
$$

So far the proof of Theorem 4 is finished.

\vskip.3in
\noindent{\bf Acknowledgements}

In April 2012, Quanhui Yang of Nanjing Normal University gave a talk
in ``ergodic prime number theorem 2012'' seminar in Morningside
Mathematical Center of the Chinese Academy of Sciences in Beijing to
introduce some results in the paper [3]. I would like to thank
Quanhui Yang for his talk which attracts my interest to this topic.
I also thank all my colleagues and friends in this seminar for
helpful discussion.
%------------------------------------------------------------------

\bigskip

\

Institute of  Mathematics, Academia Sinica, Beijing 100190, P. R.
China

E-mail: jiach@math.ac.cn
%------------------------------------------------------------------
\end{document}